\documentclass[10pt]{amsart}

\usepackage[english]{babel}
\usepackage[T1]{fontenc}
\usepackage[utf8]{inputenc}
\usepackage{amsmath}
\usepackage{amssymb}
\usepackage{amsthm}
\usepackage{color}
\usepackage[a4paper]{geometry}
\usepackage{tikz}
\usepackage{mathdots}
\usepackage[all]{xy}
\usepackage{hyperref}
\usepackage{accents}
\usepackage{lmodern}
\usepackage[multiple]{footmisc}
\usepackage{mathrsfs}
\usepackage{varwidth}

\theoremstyle{definition}
\newtheorem{definition}{Definition}[section]

\newtheorem{rmk}[definition]{Remark}
\newtheorem{quest}[definition]{Question}

\theoremstyle{plain}
\newtheorem{theorem}[definition]{Theorem}
\newtheorem{prop}[definition]{Proposition}

\newtheorem{lem}[definition]{Lemma}

\newcommand{\mb}{\mathbb}

\newcommand{\bs}{\boldsymbol}

\begin{document}

\title{A note on bounded exponential sums}

\author{Reynold Fregoli}
\address{Department of Mathematics\\ 
Royal Holloway, University of London\\ 
TW20 0EX Egham\\ 
UK}
\email{Reynold.Fregoli.2017@live.rhul.ac.uk}

\subjclass{11L03, 11L07 ; 28A78, 40C15}
\date{\today, and in revised form ....}

\dedicatory{}

\keywords{}

\begin{abstract}
Let $A\subset\mb{N}$, $\alpha\in(0,1)$, and for $x\in\mb{R}$ let $e(x):=e^{2\pi ix}$. We set
$$S_{A}(\alpha,N):=\sum_{\substack{n\in A\\n\leq N}}e(n\alpha).$$
Recently, Lambert A'Campo proposed the following question: is there an infinite non-cofinite set $A\subset\mb{N}$ such that for all $\alpha\in(0,1)$ the sum $S_{A}(\alpha,N)$ has bounded modulus as $N\to +\infty$? In this note we show that such sets do not exist. To do so, we use a theorem by Duffin and Schaeffer on complex power series. We extend our result by proving that if the sum $S_{A}(\alpha,N)$ is bounded in modulus on an arbitrarily small interval and on the set of rational points, then the set $A$ has to be either finite or cofinite. On the other hand, we show that there are infinite non-cofinite sets $A$ such that $|S_{A}(\alpha,N)|$ is bounded for all $\alpha\in E\subset (0,1)$, where $E$ has full Hausdorff dimension and $\mb{Q}\cap (0,1)\subset E$. 
\end{abstract}

\maketitle

\section{Introduction}

We denote by $\mb{N}=\{1,2,3,\dotsc\}$ the set of natural numbers and for $x\in\mb{R}$ we denote by $e(x)$ the complex number $e^{2\pi i x}$. Let $A\subset\mb{N}$, $\alpha\in(0,1)$, and $N\in\mb{N}$. We consider the sum
$$
S_{A}(\alpha,N):=\sum_{\substack{n\in A\\n\leq N}}e(n\alpha).
$$

We observe that for all $\alpha\in(0,1)$
$$\left|\sum_{n\leq N}e(n\alpha)\right|=\left|\frac{e((N+1)\alpha)-e(\alpha)}{e(\alpha)-1}\right|\leq\frac{2}{|e(\alpha)-1|}.$$
Hence, if the set $A\subset\mb{N}$ is finite or cofinite\footnote{A set $A\subset\mb{N}$ is cofinite if $A\setminus\mb{N}$ is finite.}, the sum $S_{A}(\alpha,N)$ is bounded in modulus, i.e., for each $\alpha\in(0,1)$ there exists a constant $C_{A,\alpha}>0$, only depending on the set $A$ and the real $\alpha$, such that $|S_{A}(\alpha,N)|\leq C_{A,\alpha}$ for all $N\in\mb{N}$.

Lambert A'Campo has raised the following question.

\begin{quest}[L. A'Campo]
\label{pbl:pbl}
Are there infinite non-cofinite sets $A\subset\mb{N}$ such that for each $\alpha\in(0,1)$ there exists a constant $C_{A,\alpha}>0$ for which $|S_{A}(\alpha,N)|\leq C_{A,\alpha}$ for all $N\in\mb{N}$?
\end{quest}

Question \ref{pbl:pbl} was presented by Philipp Habegger during the problem session at the "Diophantine Approximation and Transcendence" conference held in Luminy from September 10th to 14th 2018. In this note we answer Question \ref{pbl:pbl}, by showing that such sets $A$ do not exist. More generally, we prove the following.

\begin{prop}
\label{thm:mainthm}
Let $\bs{a}:=(a_{n})_{n\in\mb{N}}$ be a sequence of complex numbers taking only finitely many values, and let
$$E(\bs{a}):=\left\{\alpha\in(0,1):\sup_{N\in\mb{N}}\left|\sum_{n\leq N}a_{n}e(n\alpha)\right|<+\infty\right\}.$$
Assume that
\begin{itemize}\vspace{2mm}
\item[$i)$] the set $E(\bs{a})$ contains an open non-empty interval;\vspace{2mm}
\item[$ii)$] the set $E(\bs{a})$ contains $\mb{Q}\cap(0,1/2]$.\vspace{2mm}
\end{itemize}
Then, the sequence $\bs{a}$ is ultimately constant.
\end{prop}
\noindent An answer to Question \ref{pbl:pbl} is provided by the case $a_{n}=\chi_{A}(n)$ and $E(\bs{a})=(0,1)$.

To make things easier, we give the following definition.

\begin{definition}
Let $E\subset(0,1)$. We say that a set $A\subset\mb{N}$ has $\textup{BES}$ (bounded exponential sums) over $E$ if for each $\alpha\in E$ there exists a constant $C_{A,\alpha}>0$ such that $|S_{A}(\alpha, N)|\leq C_{A,\alpha}$ for all $N\in\mb{N}$.
\end{definition}

We note that a set $A\subset\mb{N}$ has $\textup{BES}$ over $(0,1)$ if and only if it has $\textup{BES}$ over $(0,1/2]$. Indeed, for $A\subset\mb{N}$ and $\alpha\in(0,1)$ we have
\begin{equation}
\label{eq:symmetry}
\overline{S_{A}(\alpha,N)}=\sum_{\substack{n\leq N \\ n\in A}}\overline{e(n\alpha)}=\sum_{\substack{n\leq N \\ n\in A}}e(n(-\alpha))=\sum_{\substack{n\leq N \\ n\in A}}e(n(1-\alpha))=S_{A}(1-\alpha,N),
\end{equation}
proving that the function $S_{A}(\alpha,N)$ is bounded if and only if the function $S_{A}(1-\alpha,N)$ is bounded. This shows that condition $ii)$ in Proposition \ref{thm:mainthm} is equivalent to $\mb{Q}\cap(0,1)\subset E(\bs{a})$.

Now, we analyse the two conditions appearing in Proposition \ref{thm:mainthm}. Condition $ii)$ is clearly necessary. Indeed, the series
$$\sum_{n=0}^{+\infty}e(-pn/q)e(n\alpha)=(1-e(-p/q)e(\alpha))^{-1}$$
has finitely many complex coefficients and is only unbounded at the rational $\alpha=p/q$. If we assume $a_{n}\in\{0,1\}$, condition $ii)$ can be replaced by "for each $q\geq 2$ there exists $0<p\leq q-1$ such that $(p,q)=1$ and $p/q\in E(\bs{a})$"\footnote{Note that if we do not assume $(p,q)=1$, the result no longer holds. Consider, e.g., the rational function $(z+1)/(z^{4}-1)$. This function is unbounded only at $1$, $e(1/4)$, and $e(3/4)$, and has a power series expansion whose coefficients are not ultimately constant. However, for all even $q$ we could choose $p=q/2$, so that the hypotheisis still holds.}. This is again necessary since, e.g., the set $A=\{qn\}_{n\in\mb{N}}$ has $\textup{BES}$ over $E=(0,1)\setminus\{p/q:p=1,\dotsc,q-1\}$ for all integers $q\geq 2$. On the other hand, condition $i)$ in Proposition \ref{thm:mainthm} is not strictly necessary. To see this, one can use a slightly modified version of Theorem \ref{thm:DSperiodic} in Section \ref{sec:proofmainthm} (see \cite{Helson:OnaTheorem}) which shows that the result of Proposition \ref{thm:mainthm} still holds if we remove from the interval contained in $E(\bs{a})$ a zero Lebesgue measure set. It is then natural to ask whether the presence of an interval (up to zero measure sets) in a subset $E\subset(0,1)$ is necessary to avoid the existence of an infinite  non-cofinite set $A\subset\mb{N}$ with $\textup{BES}$ over $E$. In other words, is a purely measure-theoretic condition enough? Note that there are subsets $E\subset(0,1)$ such that $\mathscr{L}(E)=1-\varepsilon$ ($0<\varepsilon<1$) and $\mathscr{L}(E\cap I)<\mathscr{L}(I)$ for any interval $I\subsetneq(0,1)$, here $\mathscr{L}$ stands for the Lebesgue measure\footnote{An example of such a set could be the following. Assume that $\mb{Q}\cap(0,1)=\{q_{n}\}_{n\geq 1}$ is a numbering of the rational numbers and let $0<\varepsilon<1$. Consider the set $E=\bigcup_{n\geq 1}(q_{n}-\varepsilon 2^{-n},q_{n}+\varepsilon 2^{-n})$. We have $\mathscr{L}((0,1)\setminus E)\geq 1-\varepsilon$. Moreover, since every non-empty interval $I\subset(0,1)$ contains a rational, we have $\mathscr{L}(((0,1)\setminus E)\cap I)<\mathscr{L}(I)$.}. In view of this, we study subsets $\mb{Q}\cap(0,1/2]\subset E\subset(0,1)$ that admit infinite non-cofinite sets $A\subset\mb{N}$ with $\textup{BES}$ over $E$. "How big" can such subsets $E$ be? A partial answer is provided by the following.

\begin{prop}
\label{thm:thm}
There exist infinite non-cofinite sets $A\subset\mb{N}$ with $\textup{BES}$ over a subset $\mb{Q}\cap(0,1/2]\subset E\subset (0,1)$ of full Hausdorff dimension. 
\end{prop} 

Proposition \ref{thm:thm} is a consequence of the following result.

\begin{prop}
\label{prop:prop1}
Let $f:\mb{N}\to\mb{N}$ be a strictly increasing function. Then, the set $A(f):=\{n+f(n)!:n\in\mb{N}\}$ has $\textup{BES}$ over $E=\mb{Q}\cap(0,1/2]$. Moreover, any function $f$ such that\vspace{2mm}
\begin{itemize}
\item[$i)$] $\sum_{i\geq 1}1/f(i)<+\infty$,\vspace{2mm}
\item[$ii)$] $\sup_{i\in\mb{N}}\left(1/i!\right)^{\varepsilon}\prod_{f(j)\leq i}(f(j)+1)<+\infty$ for all $0<\varepsilon<1$,\vspace{2mm}
\end{itemize}
gives raise to a set $A(f)$ that has $\textup{BES}$ over some subset $\mb{Q}\cap(0,1/2]\subset E(f)\subset (0,1)$ of full Hausdorff dimension.
\end{prop}
\noindent A function that satisfies both $i)$ and $ii)$ is $f(n)=n^{2}$. We give more details in Section \ref{sec:example}.

In view of the above discussion, we are led to the following questions.
\begin{quest}
\label{quest:q4}
\ 
\begin{itemize}
\item[$a)$] Are there any positive Lebesgue measure subsets $\mb{Q}\cap(0,1/2]\subset E\subset(0,1)$ that admit infinite non-cofinite sets $A\subset\mb{N}$ with $\textup{BES}$ over $E$?
\item[$b)$] Are there any zero Lebesgue measure subsets $\mb{Q}\cap(0,1/2]\subset E\subset(0,1)$ that admit no infinite non-cofinite sets $A\subset\mb{N}$ with $\textup{BES}$ over $E$?
\end{itemize}
\end{quest}
\noindent The techniques used in this note do not seem powerful enough to tackle Question \ref{quest:q4}.

A closely related question to Question \ref{pbl:pbl} was studied by Lesigne and Petersen in \cite{LesignePetersen:Bounded}, where they prove the following result.

\begin{theorem}[Lesigne-Petersen]
\label{thm:LP}
There are no sequences $\bs{a}:=\{a_{k}\}_{k\in\mb{Z}}$ with $a_{k}\in\{\pm 1\}$ such that
\begin{equation}
\label{eq:LPpbl}
\sup_{\substack{m,n\in\mb{Z} \\ n\geq 0}}\left|\sum_{k=m}^{m+n}a_{k}e^{-ik\theta}\right|\leq c(\theta)
\end{equation}
for all $\theta\in[-\pi,\pi)$ ($c(\theta)$ being a positive real constant depending on $\theta$).
\end{theorem}

To prove Theorem \ref{thm:LP}, Lesigne and Petersen consider the compact metric space $[-1,+1]^{\mb{Z}}$ (endowed with the product distance) and the shift endomorphism $\sigma$. They fix a sequence $\bs{a}\in[-1,+1]^{\mb{Z}}$ satisfying (\ref{eq:LPpbl}), and they set $X$ to be the topological closure of the orbit of $\bs{a}$ under $\sigma$. By using the spectral theorem for Hilbert spaces, they prove that any shift-invariant probability measure $\mu$ defined on $X$ (whose existence is guaranteed by the Bogolyubov-Krylov Theorem \cite[Theorem 1.1]{Sinai:DynSistII}) must be concentrated on the point $\bs{0}$, i.e., the sequence given by all zeroes. This is clearly never true when $\bs{a}\in\{\pm 1\}^{\mb{Z}}$.

We note that the hypothesis in Theorem \ref{thm:LP} is slightly different from that of Question \ref{pbl:pbl}, the key difference being the fact that the sum in (\ref{eq:LPpbl}) is bounded also for $\theta=0$. After carefully reading Lesigne and Petersen's proof, we believe that their argument can be applied to show that, once the constraint for $\theta=0$ is removed, any shift invariant probability measure $\mu$ defined on the closure $X$ of the $\sigma$-orbit of a sequence $\bs{a}\in\{0,1\}^{\mb{Z}}$ satisfying (\ref{eq:LPpbl}) must be concentrated either on the point $\bs{0}$ or on the point $\bs{1}$. In this case we say (using the terminology from \cite{LesignePetersen:Bounded}) that $\bs{a}$ is essentially $0$ or essentially $1$. This, however, does not imply that the sequence $\bs{a}$ is eventually constant. Indeed, it is easy to see that the set $A=\{n+n!:n\in\mb{N}\}$ has an essentially zero indicator function $\chi_{A}:\mb{Z}\to\{0,1\}$. To show this, it is sufficient to observe that the closure of the orbit of $\chi_{A}$ under $\sigma$ is the set
$$X=\left\{\bs{0},\sigma^{n}\left(\chi_{A}\right),\sigma^{n}\left(\chi_{\{0\}}\right):n\in\mb{Z}\right\}.$$

\noindent We conclude the introduction by noting that Proposition \ref{thm:mainthm} provides a simpler and more elementary proof of Theorem \ref{thm:LP}. 

\section{Proof of Proposition \ref{thm:mainthm}}
\label{sec:proofmainthm}

To prove Proposition \ref{thm:mainthm} we use the following powerful result by Duffin and Schaeffer \cite[Part II, Theorem I]{DuffinSchaeffer:PowerSeries}.

\begin{theorem}[Duffin-Schaeffer]
\label{thm:DSperiodic}
Let
$$u(z):=\sum_{n=0}^{+\infty}b_{n}z^{n}$$
be a power series defined over the unit disk $D:=\{|z|<1\}\subset\mb{C}$. Assume that the coefficients $b_{n}\in\mb{C}$ take only finitely many different values. Then, if the series $u(z)$ is bounded in a sector $S:=\{\theta_{1}\leq\textup{arg}(z)\leq\theta_{2},\ |z|<1\}$ of the disk $D$, where $0\leq\theta_{1}<\theta_{2}\leq 2\pi$, the sequence $\{b_{n}\}$ is ultimately periodic\footnote{This is not explicitly stated in the theorem, but it is stated at the end of the proof (see \cite[Part II, Section 4]{DuffinSchaeffer:PowerSeries}).}.
\end{theorem} 

Let $I$ be an open interval contained in $E(\bs{a})$. We consider the function $f:I\to[0,+\infty)$ defined by
$$f(\alpha):=\sup_{N\in\mb{N}}\left|\sum_{n\leq N}a_{n}e(n\alpha)\right|.$$
This is a Baire class 1 function since it is the point-wise limit\footnote{Note that the supremum of a sequence of continuous functions $\{f_{m}\}$ can be turned into a limit by considering the continuous functions $f_{M}:=\sup_{m\leq M}f_{m}$.} of a sequence of continuous functions (see \cite[Definition 11.1]{vanRooijSchikhof:ASecondCourse}). 

By \cite[Theorem 11.4]{vanRooijSchikhof:ASecondCourse}, we know that the set of continuity points of such functions is dense in their domain.
Hence, $f$ has a continuity point $P$ in $I$. This means that we can find an interval $(\alpha_{1},\alpha_{2})$ around $P$ such that the image $f((\alpha_{1},\alpha_{2}))$ is contained in a small interval around $f(P)$. It follows that $f$ is bounded in $(\alpha_{1},\alpha_{2})$ by some constant $M>0$.

For $z\in D:=\{|z|<1\}$ we let
$$u(z):=\sum_{n=0}^{+\infty}a_{n}z^{n},$$
where $a_{0}:=0$. By applying Abel's summation formula, we find that for all $\alpha\in(\alpha_{1},\alpha_{2})$, all $0\leq r<1$, and all integers $A\geq 1$ it holds 
\begin{multline}
\label{eq:Abelmainproof}
\left|\sum_{n=0}^{A}a_{n}r^{n}e(n\alpha)\right|\leq\left|\sum_{n=0}^{A}a_{n}e(n\alpha)\right|r^{A}+\sum_{n=0}^{A-1}\left|\sum_{j=0}^{n}a_{j}e(j\alpha)\right|\left(r^{n}-r^{n+1}\right) \\
\leq f(\alpha)r^{A}+f(\alpha)\left(1-r^{A}\right)=f(\alpha)\leq M.
\end{multline}
Taking the limit for $A\to +\infty$, we obtain
$$\left|u(z)\right|\leq M$$
for $z\in S:=\{2\pi\alpha_{1}\leq\textup{arg}(z)\leq 2\pi\alpha_{2},\ |z|<1\}$. Hence, by Theorem \ref{thm:DSperiodic}, the sequence $\bs{a}$ is ultimately periodic.

To conclude the proof, we show that if $f(\alpha)<+\infty$ for all rational numbers $\alpha\in\mb{Q}\cap(0,1)$ (or equivalently in $\mb{Q}\cap(0,1/2]$ by (\ref{eq:symmetry})), the period of the sequence $\{a_{n}\}$ is $1$. Suppose that ultimately $\bs{a}$ has a period of length $q\geq 1$, i.e., $a_{n}=a_{n+q}$ for all $n\geq K$, where $K$ is some large integer. Then, for $z\in D$ we have
$$
u(z)=\sum_{n=0}^{K-1}a_{n}z^{n}+\sum_{n=0}^{+\infty}z^{qn+K}\left(\sum_{j=0}^{q-1}a_{K+j}z^{j}\right)=\sum_{n=0}^{K-1}a_{n}z^{n}+z^{K}\left(\sum_{j=0}^{q-1}a_{K+j}z^{j}\right)\frac{1}{1-z^{q}}.
$$
Since $|u(re^{i\alpha})|\leq f(\alpha)$ for all $0\leq r<1$ and all $\alpha\in E(\bs{a})$ (to see this, use (\ref{eq:Abelmainproof})), the function $u(z)$ cannot have a pole at a non trivial root of unity. Hence, the polynomial
$1+z+\dotsb +z^{q-1}$ must divide $\sum_{j=0}^{q-1}a_{K+j}z^{j}$, thus showing that $a_{K+j}=a_{K+j'}$ for all $j\neq j'$.

\section{Proof of Proposition \ref{prop:prop1}}
\label{sec:Proof2}

Let $f:\mb{N}\to\mb{N}$ be a strictly increasing function and let $A(f)=\{n+f(n)!:n\in\mb{N}\}$. Clearly, $A(f)$ is neither finite nor cofinite. For $N\geq 0$ we estimate the sum
\begin{equation}
\sum_{n\leq N}e((n+f(n)!)\alpha).\nonumber
\end{equation}
First, we show that this sum is bounded for all $\alpha\in\mb{Q}\cap(0,1)$. Let $\alpha:=p/q$, with $p,q\in\mb{N}$ and $q\geq 2$. Then, for all $n\geq q$ we have
$$n+f(n)!\equiv n \pmod{q}.$$
It follows that for $N\geq q$
\begin{multline}
\label{eq:sum}
\left|\sum_{n\leq N}e((n+f(n)!)\alpha)\right|\leq\left|\sum_{n< q}e((n+f(n)!)\alpha)\right|+\left|\sum_{q\leq n\leq N}e(n\alpha)\right| \\
\leq\left|\sum_{n< q}e((n+f(n)!)\alpha)\right|+\left|\sum_{n\leq q}e(n\alpha)\right|+\left|\sum_{n\leq N}e(n\alpha)\right|,
\end{multline}
and the right-hand side in (\ref{eq:sum}) is bounded for $N\to +\infty$. To prove the second part of Proposition \ref{prop:prop1}, we need the following auxilary result (see \cite[Section 2]{Weisz:Realnumbers}).

\begin{lem}
\label{prop:factoradic}
Let $\alpha\in[0,1)$. Then, there exists a sequence of integers $(s_{n}(\alpha))_{n\in\mb{N}}$ such that $0\leq s_{n}(\alpha)\leq n-1$ and
$$\alpha=\sum_{n\geq 1}\frac{s_{n}(\alpha)}{n!}.$$
The sequence $(s_{n}(\alpha))_{n\in\mb{N}}$ associated to $\alpha$ is unique, if we exclude all those sequences $s_{n}$ such that $s_{n}=n-1$ for all sufficiently large $n$. Under this limitation, the sequence $s_{n}(\alpha)$ is eventually null if and only if $\alpha\in\mb{Q}$.
\end{lem}

\begin{rmk}
\label{rmk:rmk1}
Let $N$ be a fixed integer and let $(s_{n})_{n> N}$ be a sequence of integers such that $0\leq s_{n}\leq n-1$ for all $n> N$. Then, we have
$$\sum_{n> N}\frac{s_{n}}{n!}\leq\frac{1}{N!}.$$
This follows from the equality
$$\sum_{n=N+1}^{M}\frac{n-1}{n!}=\frac{1}{N!}-\frac{1}{M!},$$
valid for all $M\geq N+1$. This important fact will be used later on in the proof.
\end{rmk}

For a real number $\alpha\in[0,1)$ we call the unique sequence $(s_{n}(\alpha))_{n\in\mb{N}}$ given by Lemma \ref{prop:factoradic} (that does not eventually coincide with $n-1$) the factoradic representation of $\alpha$ and we call the integer $s_{n}(\alpha)$ of such sequence the $n$-th factoradic digit of $\alpha$.

Let $\bs{a}=(a_{n})_{n\in\mb{N}}$ be another sequence of strictly positive integers and assume that
\begin{equation}
\label{eq:aconv}
\sum_{n\geq 1}1/a_{n}< +\infty.
\end{equation}
We consider the set
$$E(f,\bs{a}):=\left\{\alpha\in(0,1):s_{f(i)+1}(\alpha)\leq\frac{f(i)+1}{a_{i}}\ \mbox{for all }i\geq 1\right\}.$$
Note that $E(f,\bs{a})\neq\emptyset$ whenever $f$ is not the identity function. We shall show that for any function $f$ satisfying condition $i)$ and any sequence $\bs{a}$ satisfying (\ref{eq:aconv}) the set $A(f)$ has $\textup{BES}$ over $E(f,\bs{a})$, thereby proving that $A(f)$ has $\textup{BES}$ over $E=E(f,\bs{a})\cup(\mb{Q}\cap(0,1))$.

By Abel's summation formula, we have
\begin{multline}
\label{eq:Abel}
\left|\sum_{n\leq N}e((n+f(n)!)\alpha)\right| \\
\leq\left|\sum_{n\leq N}e(n\alpha)\right|\left|e(f(N)!\alpha)\right|+\sum_{n\leq N-1}\left|\sum_{i\leq n}e(i\alpha)\right|\left|e(f(n)!\alpha)-e(f(n+1)!\alpha)\right|.
\end{multline}
Hence, to bound the left-hand side of (\ref{eq:Abel}) it is enough to bound the sum
\begin{equation}
\label{eq:difference}
\sum_{n\leq N-1}\left|e(f(n)!\alpha)-e(f(n+1)!\alpha)\right|\leq 2\sum_{n\leq N}\left|e(f(n)!\alpha)-1\right|.
\end{equation}
Let $\{\theta\}$ denote the fractional part of any real number $\theta >0$. By using the inequality $|e(\theta)-1|\leq 2\pi\{\theta\}$ (valid for $\theta\in[0,+\infty)$), we obtain
\begin{equation}
\label{eq:pi}
\sum_{n\leq N}\left|e(f(n)!\alpha)-1\right|\leq 2\pi\sum_{n\leq N}\{f(n)!\alpha\}.
\end{equation}
Now, since $\alpha\in E(f,\bs{a})$ and $f(n)\geq n$, we have
\begin{align}
\label{eq:sums}
\{f(n)!\alpha\} & =\left\{f(n)!\sum_{i\geq 1}\frac{s_{i}(\alpha)}{i!}\right\}
=\left\{\frac{s_{f(n)+1}(\alpha)}{f(n)+1}+\frac{s_{f(n)+2}(\alpha)}{(f(n)+1)(f(n)+2)}+\dotsb\right\} \\
 & \leq\frac{s_{f(n)+1}(\alpha)}{f(n)+1}+\frac{s_{f(n)+2}(\alpha)}{(f(n)+1)(f(n)+2)}+\dotsb\nonumber \\
 & \leq\frac{1}{a_{n}}+\frac{1}{f(n)+1}+\frac{1}{(f(n)+1)(f(n)+2)}+\dotsb\nonumber \\
 & \leq\frac{1}{a_{n}}+\frac{e}{f(n)+1}.\nonumber
\end{align}
Thus, combining (\ref{eq:Abel}),(\ref{eq:difference}),(\ref{eq:pi}), and (\ref{eq:sums}), we get
$$\left|\sum_{n\leq N}e((n+f(n)!)\alpha)\right|\leq\frac{2}{|e(\alpha)-1|}\left(1+4\pi\sum_{n\leq N}\left(\frac{1}{a_{n}}+\frac{e}{f(n)+1}\right)\right).$$
By (\ref{eq:aconv}) and condition $i)$, the right hand side is bounded, proving the claim.

Now, we show that the set $E(f,\bs{a})\cup(\mb{Q}\cap(0,1))$ has full Hausdorff dimension whenever the function $f$ satisfies condition $ii)$. To give a lower bound for the Hausdorff dimension of $E(f,\bs{a})\cup(\mb{Q}\cap(0,1))$ we use the so called mass distribution principle (see \cite[Principle 4.2]{Falconer:Fractal}).

\begin{lem}
\label{lem:massdistr}
Let $\mu$ be a probability measure supported on a bounded subset $X$ of $\mb{R}$. Suppose that there are strictly positive constants $a$, $s$ and $\ell_{0}$ such that
\begin{equation}
\label{eq:massdistr}
\mu(B)\leq a|B|^{s}
\end{equation}
for any interval\footnote{Note that it is enough to consider intervals since the ball-defined Hausdorff dimension coincides with the classical Hausdorff dimension (see \cite[Section 2.4]{Falconer:Fractal})} $B$ of length $|B|\leq\ell_{0}$. Then, $\textup{dim}(X)\geq s$, where $\textup{dim}$ denotes the Hausdorff dimension of a set. 
\end{lem}

We can take $X$ to be $E(f,\bs{a})\cup(\mb{Q}\cap[0,1])$, since adding a finite number of points to a set does not change its Hausdorff dimension. To apply Lemma \ref{lem:massdistr}, we need to construct a probability measure $\mu$ whose support is contained in $X$. We use a standard limit procedure to define $\mu$ (see \cite[Proposition 1.7]{Falconer:Fractal}). For $i\in\mb{N}$ we let $\rho_{i}:=1/(i!)$ and
$$Z_{i}:=\left\{\alpha\in[0,1):s_{j}(\alpha)=0\ \mbox{for }j>i\right\}.$$
First, we observe that
\begin{equation}
\label{eq:*}
E(f,\bs{a})\cup(\mb{Q}\cap[0,1])\supset\bigcap_{i\in\mb{N}}\bigcup_{\alpha\in (E(f,\bs{a})\cup\{0\})\cap Z_{i}}[\alpha,\alpha+\rho_{i}].
\end{equation}
Indeed, by definition, for each $\alpha$ lying in the right-hand side of (\ref{eq:*}) and each $i\in\mb{N}$ there exists $\alpha_{i}\in(E(f,\bs{a})\cup\{0\})\cap Z_{i}$ such that $\alpha\in[\alpha_{i},\alpha_{i}+\rho_{i}]$. This means that either $\alpha=\alpha_{i}+\rho_{i}\in\mb{Q}\cap[0,1]$ (by Remark \ref{rmk:rmk1}) or  $0\leq\alpha-\alpha_{i}<\rho_{i}$, i.e., $\alpha-\alpha_{i}$ is a number between $0$ and $1$ such that $s_{j}(\alpha-\alpha_{i})=0$ for $j\leq i$. Thus, when we add $\alpha-\alpha_{i}$ to $\alpha_{i}$ the digits before the $i$-th do not change, showing that $s_{j}(\alpha)=s_{j}(\alpha_{i})$ for $j\leq i$. It follows that $\alpha\in E(f,\bs{a})\cup\{0\}$.

Now, for all $i\in\mb{N}$ and all $\alpha\in (E(f,\bs{a})\cup\{0\})\cap Z_{i}$ we define
\begin{equation}
\label{eq:mudef}
\mu((\alpha,\alpha+\rho_{i})):=\frac{1}{\#((E(f,\bs{a})\cup\{0\})\cap Z_{i})},
\end{equation}
where we take open intervals to make sure that for any fixed $i$ all the sets $(\alpha,\alpha+\rho_{i})$ are disjoint.
\begin{rmk}
Note that for $\alpha\in(E(f,\bs{a})\cup\{0\})\cap Z_{i}$ we have
$$\#\{\beta\in(E(f,\bs{a})\cup\{0\})\cap Z_{i+1}:\beta_{i}=\alpha\}=\frac{\#((E(f,\bs{a})\cup\{0\})\cap Z_{i+1})}{\#((E(f,\bs{a})\cup\{0\})\cap Z_{i})},$$
where $\beta_{i}$ is the truncation of $\beta$ at the $i$-th digit. Hence,
\begin{align}
\label{eq:(12)}
\mu((\alpha,\alpha+\rho_{i})) & =\frac{1}{\#((E(f,\bs{a})\cup\{0\})\cap Z_{i})}=\sum_{\substack{\beta\in(E(f,\bs{a})\cup\{0\})\cap Z_{i+1} \\ \beta_{i}=\alpha}}\frac{1}{\#((E(f,\bs{a})\cup\{0\})\cap Z_{i+1})}\nonumber \\
 & =\sum_{\substack{\beta\in(E(f,\bs{a})\cup\{0\})\cap Z_{i+1} \\ \beta_{i}=\alpha}}\mu((\beta,\beta+\rho_{i+1}))=\sum_{\substack{\beta\in(E(f,\bs{a})\cup\{0\})\cap Z_{i+1} \\ (\beta,\beta+\rho_{i+1})\subset(\alpha,\alpha+\rho_{i})}}\mu((\beta,\beta+\rho_{i+1})).
\end{align}
\end{rmk}
By \cite[Proposition 1.7]{Falconer:Fractal} and (\ref{eq:(12)}), Equation (\ref{eq:mudef}) induces a unique well-defined Borel measure $\mu$ on $\mb{R}$, with the property\footnote{Stated a few lines above \cite[Proposition 1.7]{Falconer:Fractal}} that
\begin{equation}
\label{eq:rightprop}
\mu\left(\mb{R}\setminus\bigcup_{\alpha\in (E(f,\bs{a})\cup\{0\})\cap Z_{i}}(\alpha,\alpha+\rho_{i})\right)=0
\end{equation}
for all $i\in\mb{N}$, and supported on the set
\begin{equation}
\label{eq:***}
\bigcap_{i\in\mb{N}}\bigcup_{\alpha\in (E(f,\bs{a})\cup\{0\})\cap Z_{i}}[\alpha,\alpha+\rho_{i}]\subset X.
\end{equation}

To prove (\ref{eq:massdistr}), we fix a number $i_{0}\in\mb{N}$ and a real $0<s<1$. We consider an interval $B\subset[0,1]$ of length $|B|$ less than $\rho_{i_{0}}$. Clearly, there exists an index $i\in\mb{N}$, $i\geq i_{0}$, such that
\begin{equation}
\label{eq:Bbounds}
\rho_{i+1}<|B|\leq \rho_{i}.
\end{equation}
By (\ref{eq:rightprop}), we have
$$\mu(B)\leq\mu\left(\bigcup_{\substack{\alpha\in (E(f,\bs{a})\cup\{0\})\cap Z_{i+1} \\ \left(\alpha,\alpha+\rho_{i+1}\right)\cap B\neq\emptyset}}\left(\alpha,\alpha+\rho_{i+1}\right)\right),$$
and it is straightforward to see that
$$\#\{\alpha\in Z_{i+1}:(\alpha,\alpha+\rho_{i+1})\cap B\neq\emptyset\}\leq\frac{|B|}{1/(i+1)!}+2,$$
since the intervals $(\alpha,\alpha+\rho_{i+1})$ are pairwise disjoint and each of them has length $1/(i+1)!$. Hence, we have
\begin{multline}
\label{eq:Hausdim}
\mu(B)\leq \sum_{\substack{\alpha\in (E(f,\bs{a})\cup\{0\})\cap Z_{i+1} \\ (\alpha,\alpha+\rho_{i+1})\cap B\neq\emptyset}}\mu((\alpha,\alpha+\rho_{i+1})) \\
\leq\left(\frac{|B|}{1/(i+1)!}+2\right)\frac{1}{\#((E(f,\bs{a})\cup\{0\})\cap Z_{i+1})}.
\end{multline}
Now, we observe that for all $j\in\mb{N}$
\begin{equation}
\label{eq:Hausdorffmeasure1}
\#((E(f,\bs{a})\cup\{0\})\cap Z_{j})\geq\frac{j!}{\prod_{f(k)+1\leq j}(f(k)+1)},
\end{equation}
since for $\alpha\in E(f,\bs{a})\cup\{0\}$ the $0$ digit is always allowed in the $(f(k)+1)$-th position independently of $\bs{a}$. Thus, by (\ref{eq:Bbounds}), (\ref{eq:Hausdim}), and (\ref{eq:Hausdorffmeasure1}), we deduce
\begin{multline}
\mu(B)\leq\left(\frac{|B|}{1/(i+1)!}+2\right)\frac{\prod_{f(j)+1\leq i+1}(f(j)+1)}{(i+1)!}\leq 3(i+1)!|B|\frac{\prod_{f(j)\leq i}(f(j)+1)}{(i+1)!} \\
\leq 3|B|^{1-s}\prod_{f(j)\leq i}(f(j)+1)|B|^{s}\leq 3\left(\frac{1}{i!}\right)^{1-s}\prod_{f(j)\leq i}(f(j)+1)|B|^{s}.\nonumber
\end{multline}
It then follows from Lemma \ref{lem:massdistr} that any function $f$ such that
$$\sup_{i\in\mb{N}}\left(\frac{1}{i!}\right)^{1-s}\prod_{f(j)\leq i}(f(j)+1)<+\infty$$
for all $0<s<1$ gives raise to a set $E(f,\bs{a})$ of full Hausdorff dimension.

\section{An example}
\label{sec:example}
To conclude this note, we give an example of a function $f:\mb{N}\to\mb{N}$ satisfying both conditions $i)$ and $ii)$ of Proposition \ref{prop:prop1}. We let $s:=1-\varepsilon$ and $f(i)=i^{2}$ for $i\in\mb{N}$. Condition $i)$ is clearly satisfied. Moreover, we have
\begin{equation}
\label{eq:Hausdorffmeasure}
\left(\frac{1}{i!}\right)^{1-s}\prod_{f(j)\leq i}(f(j)+1)=\frac{\prod_{j\leq \sqrt{i}}(j^{2}+1)}{(i!)^{\varepsilon}}.
\end{equation}
Now, when $i>\lceil 3/\varepsilon\rceil\left\lfloor\sqrt{i}\right\rfloor$, we find
$$(i!)^{\varepsilon}\geq\left(\underbrace{1\dotsm 1}_{\lceil 3/\varepsilon\rceil\ \mbox{times}}\cdot\underbrace{2\dotsm 2}_{\lceil3/\varepsilon\rceil\ \mbox{times}}\cdot\dotsm\cdot\underbrace{\left\lfloor\sqrt{i}\right\rfloor\dotsm\left\lfloor\sqrt{i}\right\rfloor}_{\lceil3/\varepsilon\rceil\ \mbox{times}}\right)^{\varepsilon}\geq\prod_{j\leq\sqrt{i}}j^{3}.$$
Hence, the right-hand side in (\ref{eq:Hausdorffmeasure}) is always bounded, showing that $ii)$ holds.

\section*{Acknowledgements}
I am extremely grateful to Jeffrey Vaaler for pointing me to Helson's paper \cite{Helson:OnaTheorem} and for the very fruitful discussions that we had during his visit at Royal Holloway. I would like thank my supervisor, Martin Widmer, who attended the conference in Luminy and introduced me to this problem, for his encouragement and precious advice. I would also like to thank Lambert A'Campo and Cedric Pilatte for spotting a mistake in an early version of this manuscript. Finally, I am thankful to Royal Holloway University of London, for funding my position here.

\addcontentsline{toc}{section}{\bibname}
\bibliographystyle{plain}
\bibliography{Bibliography}

\end{document}